\documentclass[11pt,leqno]{article}
\usepackage{amssymb}
\usepackage{amsmath}

\newtheorem{theorem}{Theorem}[section]

\numberwithin{equation}{section}

\begin{document}

\title{Time Decay for solutions to One Dimensional Two Component Plasma Equations}
\author{Robert Glassey \\
Stephen Pankavich \\
Jack Schaeffer}
\date{\today}
\maketitle

\begin{center}
\end{center}


\section{Dedication and Introduction}
We represent three generations of students: Bob Glassey, Walter's student finishing at Brown in
1972; Jack Schaeffer, Bob's student finishing at Indiana University in 1983 and Steve Pankavich,
Jack's student finishing at Carnegie Mellon in 2005.  We have all thrived professionally from
our association with Walter and are delighted to dedicate this note to him on the occasion of
his 70th birthday. The problem we study below concerns the asymptotic behavior of solutions, an
area to which Walter has contributed greatly.

The motion of a collisionless plasma is described by the
Vlasov--Maxwell system.  If we neglect magnetic effects we then have the
Vlasov--Poisson system (VP).  We can also consider the effect of large velocities
and solutions to the relativistic Vlasov--Poisson system (RVP).
We will study both systems in one space and one momentum dimension,
with two species of oppositely charged particles.  We further assume that
each system is {\it neutral\/} which means that the average
value of the density $\rho$ vanishes (see below).
The Vlasov--Poisson system (VP) is
\begin{equation}\label{VP}
\left \{ \begin{gathered}
\partial_t f + v \ \partial_x f  + E(t,x) \ \partial_v f = 0,\\
\partial_t g + \frac{v}{m} \ \partial_x g  - E(t,x) \ \partial_v g = 0,\\
\rho(t,x) = \int \left ( f(t,x,v) - g(t,x,v) \right ) \ dv, \\
E(t,x) = \frac{1}{2} \left ( \int_{-\infty}^x \rho(t,y) \ dy -
\int_x^\infty \rho(t,y) \ dy \right ).
\end{gathered} \right.
\end{equation}

Here $t \geq 0$ is time, $x \in \mathbb{R}$ is position, $v \in
\mathbb{R}$ is momentum, $f$ is the number density in phase space of
particles with mass one and positive unit charge while  $g$ is the
number density of particles with mass $m > 0$ and negative unit
charge. The effect of collisions is neglected.  The initial
conditions $$ f(0,x,v) = f_0(x,v) \geq 0$$ and $$g(0,x,v) = g_0(x,v)
\geq 0$$ for $(x,v) \in \mathbb{R}^2$ are prescribed.  We assume
that $f_0,\,g_0 \in C^1(\mathbb{R}^2)$ are nonnegative, compactly
supported and satisfy the neutrality condition
\begin{equation}
\label{neutrality} \iint f_0 \ dv dx = \iint g_0 \ dv\, dx.
\end{equation}
Using the notation
$$ \hat{v}_m = \frac{v}{\sqrt{m^2 + v^2}},$$ we can write the relativistic
Vlasov--Poisson system (abbreviated RVP) as
\begin{equation}
\label{RVP} \left \{ \begin{gathered}
\partial_t f + \hat{v}_1 \ \partial_x f  + E \ \partial_v f = 0,\\
\partial_t g + \hat{v}_m \ \partial_x g  - E \ \partial_v g = 0,\\
\rho(t,x) = \int \left ( f - g \right ) \ dv, \\
E(t,x) = \frac{1}{2} \left ( \int_{-\infty}^x \rho \ dy -
\int_x^\infty \rho \ dy \right ).
\end{gathered} \right.
\end{equation}

Global existence and regularity are known for solutions of (\ref{VP}) and (\ref{RVP}).  Both
$f(t,\cdot, \cdot)$ and $g(t,\cdot,\cdot)$ are compactly supported for all $t \geq 0$. There is
scant literature regarding the large time behavior of solutions.  Some time decay is known for
the three-dimensional analogue of (\ref{VP})(\cite{GS}, \cite{IR}, \cite{Per}).  Also, there are
time decay results for (\ref{VP}) (in dimension one) when the plasma is monocharged (set $g
\equiv 0$) (\cite{BKR}, \cite{BFFM}, \cite{Sch}).  In this work two species of particles with
opposite charge are considered, thus the methods used in these references do not apply.
References \cite{DD}, \cite{Dol}, and \cite{DR} are also mentioned since they deal with
time-dependent rescalings and time decay for other kinetic equations.  We will take $m=1$ below.
A full description of these results will appear in \cite{GPS}.

First we sketch the derivation of an identity for solutions to (\ref{VP}) from which we can
conclude that certain positive quantities are integrable in $t$ on the interval $[0,\infty)$.
This identity also extends to (\ref{RVP}) but the results are weaker.  Unfortunately, these
identities are very ``one-dimensional'', that is, they do not seem to easily generalize to
higher dimension. Moreover, it is not clear if there is an extension which allows for more than
two species of particles.


Here are the results we have obtained.  The classical equations for (VP) are
$$f_t+vf_x+Ef_v=0,\quad g_t+vg_x-Eg_v=0$$
where $E_x=\rho=\int (f-g)\,dv$.
Let
$$F(t,x)=\int f(t,x,v)\,dv, \quad  G(t,x)=\int g(t,x,v)\,dv.$$
Then $\rho=F-G$.  We will show that
$$\int_0^{\infty}\!\!\int E^2 (F+G)\,dx\,dt<\infty.$$
From this it will follow in the non--relativistic case that
$$\int_0^{\infty}\!\!\int (F^4+G^4)\,dx\,dt<\infty$$
while the corresponding result for solutions to (RVP) is
$$\int_0^\infty \left (\int \left(F(t,x)^\frac{7}{4}+ G(t,x)^\frac{7}{4}\right)\ dx
\right )^4 \ dt < \infty.$$ The local charges for solutions to both systems will satisfy for any
fixed $R>0$
$$\lim_{t\to \infty}\int_{|x|<R}F(t,x)\,dx=
\lim_{t\to \infty}\int_{|x|<R}G(t,x)\,dx=0.$$ Finally, for solutions to \ref{VP} or \ref{RVP} we
can show that
$$\lim_{t\to \infty}\|E(t,\cdot)\|_{\infty}=0.$$

\section{Results}
We first derive a general identity which holds for both (VP) and (RVP).  From the above
definitions and (VP) we have
$$F_t=-\int (vf_x+Ef_v)\,dv = -\partial_x \int vf\,dv$$
and thus
$$\partial_t \int_{-\infty}^xF(t,y)\,dy = -\int vf(t,x,v)\,dv$$
with a similar result for $g$.
Multiply the $f$ equation in  (VP) by $v\cdot \int_{-\infty}^xF(t,y)\,dy$
and integrate over $v$:
$$\int vf_t \int_{-\infty}^xF(t,y)\,dy \,dv +
\int v^2f_x \int_{-\infty}^xF(t,y)\,dy \,dv +
\int vf_vE \int_{-\infty}^xF(t,y)\,dy \,dv=0.$$
Write this as $I+II+III=0$.  Then
\begin{eqnarray}
I&=&\partial_t \left[\int vf \int_{-\infty}^xF(t,y)\,dy \,dv\right]-
\int vf \,dv\int_{-\infty}^xF_t(t,y)\,dy\\
&=&\partial_t \left[\int vf \int_{-\infty}^xF(t,y)\,dy \,dv\right ]+ \left(\int vf\,dv\right)^2,
\end{eqnarray}
$$II=\partial_x \left[\int v^2f \int_{-\infty}^xF(t,y)\,dy \,dv\right]-
\int v^2f\cdot F(t,x)\,dv,$$ and after integrating by parts in $v$
$$III=-\partial_x \left[\frac12 E(t,x)\left(\int_{-\infty}^xF(t,y)\,dy\right)^2\right]
+\frac12\rho(t,x)\left[\int_{-\infty}^xF(t,y)\,dy\right]^2.$$
Now integrate over $x$:
$$\frac{d}{dt} \int\left[\int vf \int_{-\infty}^xF(t,y)\,dy \,dv\right]\,dx
+\int \left(\int vf\,dv\right)^2\,dx\\ -\int F(t,x)\int v^2f\,dv\,dx$$
$$+\frac12 \int \rho(t,x)\left[\int_{-\infty}^xF(t,y)\,dy\right]^2\,dx=0.$$

Now repeat this calculation with $f$ replaced by $g$ and add the two
results to derive
$$\frac{d}{dt} \int\left[\int vf \int_{-\infty}^xF(t,y)\,dy \,dv+
\int vg \int_{-\infty}^xG(t,y)\,dy \,dv\right]\,dx$$
$$ +\int \left(\int vf\,dv\right)^2\,dx\\ -\int F(t,x)\int v^2f\,dv\,dx$$
$$+\int \left(\int vg\,dv\right)^2\,dx\\ -\int G(t,x)\int v^2g\,dv\,dx $$
$$+\frac12\int \rho(t,x)\left(\left[\int_{-\infty}^xF(t,y)\,dy\right]^2-
\left[\int_{-\infty}^xG(t,y)\,dy\right]^2\right)dx
=0.$$
The first line is bounded when integrated in time.  The second and third lines
are nonpositive.
Call $L$ the last term above.  Then because $\rho = \int (f-g)\,dv=F-G$ and
$E=\int_{-\infty}^x\rho(t,y)\,dy=\int_{-\infty}^x(F-G)\,dy$ we get after
a brief calculation
$$L=-\frac14 \int E^2 (F+G)\,dx.$$

Thus in particular
$$\int_0^{\infty}\!\!\int E^2 (F+G)\,dx\,dt<\infty$$
and
$$\int_0^{\infty}\!\!\int \left[F(t,x)\int v^2f\,dv-
\left(\int vf\,dv\right)^2\right]\,dx\,dt<\infty$$
$$\int_0^{\infty}\!\!\int \left[G(t,x)\int v^2g\,dv-
 \left(\int vg\,dv\right)^2\right]\,dx\,dt<\infty.$$

We can use these inequalities directly to establish the $L^4$ estimate. Write
$$F(t,x)\int v^2f\,dv- \left(\int vf\,dv\right)^2$$
as
$${1\over 2}\int \int  (w-v)^2f(v) f(w)\,dv\,dw.$$
Then, from above we know that the quantity
$$k(t,x)\equiv \int \int  (w-v)^2f(t,x,v) f(t,x,w)\,dv\,dw$$
is integrable over all $x,\,t$.  To get the $L^4$ bound we split the integral
for $F(t,x)^2$ in the usual manner:
$$F(t,x)^2=\int \int f(v) f(w)\,dv\,dw=\int_{|v-w|<R}+\int_{|v-w|>R}\equiv
I_1+I_2.$$
So $I_2\le R^{-2}k(t,x)$ and in $I_1$
$$\int_{|v-w|<R}f(w)\,dw=\int_{v-R}^{v+R}f(w)\,dw\le cR.$$
Thus
$$I_1\le c\cdot R \cdot F.$$
Set $R=k^{1/3}F^{-1/3}$.  Then
$F^4\le ck$ so $F^4$ is integrable over all $x,\,t$.
The result for $G$ is exactly the same.

Our final results will show that for solutions to the classical VP
system (\ref{VP}) and RVP system (\ref{RVP}), the electric field $E$
tends to 0 in the maximum norm.  Here we consider the former only.

\begin{theorem} \label{E3int}
Under the above assumptions  consider
solutions $f,\,g$ to (\ref{VP}). Then
$$\int_0^{\infty}\Vert E(t) \Vert_{\infty}^3\,dt<\infty.$$
\end{theorem}
{\bf Proof:}  This will follow immediately from the  result
 that
$$Q(t) := \int_{-\infty}^{\infty} E^2(t,x)\Big[F(t,x)+G(t,x)\Big]\,dx$$
is integrable in time.  Indeed because $E_x=\rho =
\int(f-g)\,dv=F-G$, we have
$$\frac{\partial}{\partial x}E^3=3E^2\rho=3E^2(F-G).$$
Integrate in $x$ to get
$$E^3(t,x)=\int_{-\infty}^x 3E^2(F-G)\,dx$$
so that
\begin{equation}
\label{E3int} |E(t,x)|^3\le \int_{-\infty}^{\infty} 3E^2(F+G)\,dx
= 3Q(t)
\end{equation}
and the result follows as claimed.  This can now be exploited to show
that the electric field tends uniformly to 0 as $t\to \infty$.

\begin{theorem} \label{EdecayVP}
Let the previous assumptions hold and consider
solutions $f,\,g$ to the classical (VP) system (\ref{VP}).   Then
$$\lim_{t\to \infty}\|E(t)\|_{\infty}=0.$$
\end{theorem}
{\bf Proof:}  We will show that
$$\lim_{t\to \infty}Q(t)=0.$$
The conclusion will then follow from (\ref{E3int}). Since $Q(t)$
is integrable over $[0,\infty)$, $\liminf Q(t)=0$.  If we can show
that $\dot Q(t)$ is bounded we may then conclude the statement of the theorem.
(An alternate proof is given in \cite{GPS}.)

Using $E_x=\rho=F-G$ and $E_t=-j=-\int v(f-g)\,dv$ we
compute
$$
\frac{dQ}{dt}= -2\int jE(F+G)\,dx + 2\int \rho E\int v(f+g)\,dv\,dx.$$
Now, $E$ is uniformly bounded because by definition in (\ref{VP})
$$|E(t,x)|\le \int_{-\infty}^x (F+G)(t,x)\,dx\le
\int_{-\infty}^{\infty} (F+G)(t,x)\,dx\le \hbox{const.}$$ where the
last inequality follows by conservation of mass.   Therefore
$$\left|\frac{dQ}{dt}\right| \le  c \int (F+G)\int |v|(f+g)\,dv \,dx.$$
Define $e$ to be the kinetic energy density, $$e(t,x) := \int v^2
(f+g)\,dv.$$ Then by splitting the $v$ integral into sets $|v|<R$ and its
complement we get
$$\int |v|(f+g)\,dv \le  c e^\frac{2}{3}(t,x).$$
Similarly
$$F\equiv \int f \,dv \le  c e^\frac{1}{3}(t,x)$$
and therefore
$$\left|\frac{dQ}{dt}\right|\le  c \int (F+G)e^\frac{2}{3}\,dx \le
c \Big(\int (F+G)^3 \,dx\Big)^{1\over 3}\le \hbox{const.}$$ This concludes the proof.  Complete
details may be found in \cite{GPS}.

\end{document}